\documentclass[francais]{article} 
\usepackage[latin1]{inputenc}
\usepackage{fontenc}
\usepackage{amsmath}
\usepackage{amssymb}
\usepackage{graphicx}
\usepackage{epsfig}
\usepackage{babel}
 
\newtheorem{theorem}{ThÈorËme}

\def\epsilon{\varepsilon}
\def\R{{\mathbb R}} 
 
\def\Z{{\mathbb Z}}

\def\G{{\cal G}} 
\def\A{{\cal A}} 
\def\F{{\cal F}} 
 
\def\1{u} 
\def\S{\mathfrak{S}}
\def\g{\tilde{g}}

\def\s{\text{sys}}

\def\epsilon{{\varepsilon}}

\title{Sur la systole de la sphËre au voisinage de la mÈtrique standard}

\author{Florent BALACHEFF\footnote{Section de MathÈmatiques, UniversitÈ de GenËve, rue du LiËvre 2-4, C.P. 64 CH-1211 Gen\`eve 4,  Mail :  Florent.Balacheff@math.unige.ch}} 

\begin{document}

\maketitle

\begin{abstract}
Nous Ètudions l'aire systolique (dÈfini com\-me le rapport de l'aire sur le carrÈ de la systole) de la sphËre munie d'une mÈtrique riemannienne lisse comme fonction de cette mÈtrique. Cette fonction, bornÈe infÈrieurement par une quantitÈ strictement positive sur l'espace des mÈtri\-ques, admet la mÈtrique ronde standard $g_0$ pour point critique, bien que celle-ci ne corresponde pas au minimum global conjecturÈ : nous montrons que pour toute direction tangente ‡ l'espace des mÈtriques en $g_0$, il existe une variation de mÈtrique dans cette direction le long de laquelle l'aire systolique ne peut qu'augmenter.

\end{abstract}

\bigskip

\section{Introduction}

…tant donnÈe une mÈtrique riemannienne lisse $g$ sur la sphËre de dimension $2$ notÈe $S^2$, on dÈfinit la {\it systole} de $(S^2,g)$ comme la plus petite longueur d'une gÈodÈsique fermÈe non rÈduite ‡ un point. On note $\s(S^2,g)$ cette quantitÈ, qui est strictement positive, et rÈalisÈe comme la longueur d'une gÈodÈsique fermÈe.

Nous pouvons alors dÈfinir la {\it constante systolique} de $S^2$ comme la quantitÈ
$$
\S(S^2)=\inf_g  \frac{\A(S^2,g)}{(\s(S^2,g))^2} ,
$$
o˘ $\A(S^2,g)$ dÈsigne l'aire riemannienne de $(S^2,g)$ et l'infimum est pris sur l'ensemble des mÈtriques riemanniennes lisses de $S^2$. Cette quantitÈ est non nulle : ce rÈsultat a ÈtÈ dÈmontrÈ par C. Croke \cite{crok88}, puis rÈcemment amÈliorÈ par A. Nabutovsky et R. Rotman \cite{nabrot02} d'une part, S. Sabourau \cite{sabo04} d'autre part :
$$
\S(S^2) \geq \frac{1}{64}.
$$
Il est conjecturÈ que 
$$
\S(S^2)=\frac{1}{2\sqrt{3}},
$$
 dont la valeur est atteinte pour la mÈtrique singuliËre (proposÈe par E. Calabi) obtenue en recollant le long de leur bord deux triangles ÈquilatÈraux. Notons dËs maintenant que le rapport aire sur systole au carrÈ pour la mÈtrique standard de la $2$-sphËre vaut $\frac{1}{\pi}$, ce qui est strictement supÈrieur ‡ la valeur $\frac{1}{2\sqrt{3}}$.

\bigskip

Nous nous intÈressons au comportement local du rapport de l'aire sur la systole au carrÈ pour une mÈtrique $g$ au voisinage de la mÈtrique standard $g_0$, que nous nommerons {\it aire systolique} et noterons naturellement $\S(S^2,g)$. Ce rapport n'est pas continu sur l'espace des mÈtriques lisses muni de la topologie forte, mais il est continu au voisinage de la mÈtrique standard. En effet, E. Calabi et J. Cao ont montrÈ que pour une mÈtrique $g$ ‡ courbure positive sur la sphËre $S^2$, la systole coÔnicidait avec la plus petite longueur d'une gÈodÈsique fermÈe obtenue par un procÈdÈ de minimax sur l'espace des cycles unidimensionnels (voir \cite{calcao92}). Cette derniËre quantitÈ Ètant continue en la mÈtrique, et les mÈtriques proches de la mÈtrique standard Ètant ‡ courbure positive, nous obtenons la continuitÈ de l'aire systolique au voisinage de la mÈtrique standard.

\bigskip

…tant donnÈe une variation $\{g_t\}$ de $g_0$ par une famille de mÈtriques lisses dÈpendant de maniËre lisse d'un paramËtre rÈel $t$, nous Ètudions le comportement possible de $\S(S^2,g_t)$ en fonction de $t$ au voisinage de $0$. Par le thÈorËme d'uniformisation, toute mÈtrique sur $S^2$ est conformÈment Èquivalente ‡ la mÈtri\-que standard, et pour cette raison, nous ne nous intÈresserons qu'aux variations de $g_0$ par des mÈtriques dans sa classe conforme.

L'espace des mÈtriques lisses conformes ‡ $g_0$ est naturellement identifiÈ ‡ l'espace $C^\infty(S^2,\R_+^\ast)$ des fonctions lisses sur $S^2$ ‡ valeurs rÈelles strictement positives. Nous munissons cet espace de la topologie forte. Comme cet espace est localement homÈomorphe ‡ l'espace vectoriel $C^\infty(S^2,\R)$ des fonctions lisses sur $S^2$ ‡ valeurs rÈelles (muni Ègalement de la topologie forte), ce dernier peut Ítre pensÈ comme l'espace tangent en $g_0$ ‡ l'espace des mÈtriques lisses conformes ‡ $g_0$.

L'application qui ‡ $\Phi \in C^\infty(S^2,\R^\ast_+)$ associe $\S(S^2,\Phi \cdot g_0)$ est donc une application continue au voisinage de l'application constante Ègale ‡ $1$ et nous obtenons le rÈsultat suivant :

\bigskip

\noindent {\bf ThÈorËme} {\it Pour chaque $f \in C^\infty(S^2, \R)$, il existe une famille lisse ‡ un paramËtre de fonctions $\Phi_t \in C^\infty(S^2,\R_+^\ast)$ telle que $\Phi_o=1$, $\left. \frac{d \Phi_t}{dt}\right |_0=f$ et
$$
\S(S^2,\Phi_t \cdot g_0) \geq \frac{1}{\pi},
$$
avec ÈgalitÈ pour $t\neq 0$  si et seulement si $f$ est une fonction constante.
}

\bigskip

 Cet ÈnoncÈ suggËre fortement que $g_0$ soit un minimum local.
 
\medskip

Nous pouvons gÈnÈraliser la notion de point critique d'une fonction continue dÈfinie sur une variÈtÈ compacte de dimension finie d˚e ‡ M. Morse \cite{mors59} aux fonctions continues dÈfinies sur un espace topologique localement homÈomorphe ‡ un espace vectoriel : Ètant donnÈs un espace topologique $M$ modelÈ localement sur un espace vectoriel topologique $V$, et une application $F : M \rightarrow \R$ continue, un point $x_0$ de $M$ est dit {\it rÈgulier} s'il existe une carte locale $(U,\phi)$ de classe $C^0$ centrÈe en $x_0$ et une forme linÈaire non nulle $T : V \rightarrow \R$ telle que 
$$
F(x)=F(x_0)+T(\phi(x)),
$$
pour chaque $x \in U$. Dans le cas contraire, le point est dit {\it critique}.

Nous avons le corollaire suivant :

\bigskip

\noindent {\bf Corollaire} {\it La mÈtrique $g_0$ est un point critique de la fonctionnelle aire systolique.}

\bigskip

Notons qu'il est remarquable d'obtenir un point critique, qui ne corresponde pas au minimum global conjecturÈ. 

\bigskip

Cet article est organisÈ commme suit. Dans la section suivante, nous commenÁons par montrer que l'aire systolique est un minimum local strict le long de certaines variations associÈes aux ÈlÈments $f \in C^\infty(S^2,\R)$ vÈrifiant $\int_{S_2} f dv_{g_0}=0$. Un corollaire immÈdiat en est le thÈorËme annoncÈ. Dans la troisiËme section, nous prÈsentons deux rÈsultats d˚s respectivement ‡ V. Guillemin \cite{guil76} et ‡ P. Pu \cite{pu52} qui permettent de mieux apprÈhender le comportement local de l'aire systolique au voisinage de $g_0$. Enfin, nous formulons dans la derniËre section quelques remarques et questions.

\section{…tude locale de l'aire systolique et variations minimisantes}

Le thÈorËme annoncÈ dans l'introduction est une consÈquence du rÈsultat suivant :

\bigskip

\noindent {\bf Proposition} {\it Soit $f :S^2 \rightarrow \R$ une fonction lisse non nulle vÈrifiant 
$$
\int_{S^2} f dv_{g_0} =0.
$$
 On associe ‡ $f$ la variation de mÈtrique $g_t=(1+tf)^2\cdot g_0$ et on note $a>0$ le plus grand rÈel (Èventuellement infini) tel que $g_t$ soit bien dÈfinie pour $|t|<a$. Alors il existe $\alpha \in ]0,a[$ tel que pour tout rÈel t vÈrifiant  $0<|t|<\alpha$,
$$
\S(S^2,g_t) > \S(S^2,g_0)=\frac{1}{\pi}.
$$
}

\bigskip

En effet, considÈrons la variation $g_t=(1+t)\cdot g_0$, et notons $\1 \in C^\infty(S^2,\R)$ le vecteur dÈrivÈ  correspondant ($\1$ est la fonction constante sur $S^2$ Ègale ‡ $1$). Nous avons une dÈcomposition naturelle de $C^\infty(S^2,\R)$ en somme directe :
$$
C^\infty(S^2,\R) = \R \1 \oplus \{ f \in C^\infty(S^2,\R) \mid \int_{S^2} f dv_{g_0}=0\}. \eqno (2.1)
$$
Remarquons que l'espace $\{ f \in C^\infty(S^2,\R) \mid \int_{S^2} f dv_{g_0}=0\}$ peut Ítre pensÈ comme le noyau de la diffÈrentielle en $g_0$ de l'application $\A(S^2, \cdot)$ qui ‡ $\Phi \in C^\infty(S^2,\R_+^\ast)$ associe $\A(S^2, \Phi \cdot g_0)$.

Comme l'aire systolique est invariante par changement d'Èchelle (la transformation d'une mÈtrique $g$ en $\mu \cdot g$ o˘ $\mu$ est un rÈel strictement positif ne modifie pas la valeur de l'aire systolique), nous obtenons pour tout $\lambda \in \R$ et pour toute mÈtrique $g$
$$
\S(S^2,(1+\lambda t)\cdot g)=\S(S^2,g).
$$
…tant donnÈe une fonction $f \in C^\infty(S^2,\R)$, nous pouvons dÈcomposer $f$ dans la somme directe (2.1) comme 
$$
f=\lambda(f) \1+(f - \lambda(f) \1)
$$
 o˘ $\lambda(f)=\int_{S^2} f dv_{g_0}$, et en considÈrant la variation 
$$
\g_t=(1+\lambda(f) t)(1+t(f-\lambda(f)))^2 g_0,
$$
 nous obtenons le thÈorËme.

\bigskip

Pour dÈmontrer la proposition, nous allons procÈder comme suit. …tant fixÈe la fonction lisse non nulle $f : S^2 \rightarrow \R$ vÈrifiant $\int_{S^2} f dv_{g_0} =0$, nous pouvons calculer l'aire de la mÈtrique $g_t$ associÈe. Nous voyons facilement que 
$$
\A(S^2,g_t)=\A(S^2,g_0)+t^2 \int_{S^2} f^2 dv_{g_0}.
$$
La proposition va alors dÈcouler de l'inÈgalitÈ
$$
\s(S^2,g_t) \leq 2\pi \eqno(2.2)
$$
pour $t$ dans un voisinage ouvert de $0$. Avant de dÈmontrer cette inÈgalitÈ, nous commenÁons tout d'abord par quelques rappels.

\subsection{Rappels et notations}

\noindent  {\bf Le principe du minimax.} Un outil important pour Ètudier la systole est fourni par le proc\'ed\'e de minimax. Pour plus de dÈtails, voir \cite{KLIN}.

\medskip

Soit $(S^2,g)$ une 2-sphËre riemannienne. Nous commenÁons par considÈrer $\Lambda S^2$ l'{\it espace des courbes fermÈes}, dÈfini comme l'espace des applications \linebreak $c :S^1 \rightarrow S^2$ de classe $C^1$ par morceaux muni de la topologie $C^0$. Les groupes $SO(2)$ et $O(2)$ agissent canoniquement sur $\Lambda S^2$ et on note $\Pi S^2$ et $\tilde{\Pi}S^2$ les quotients respectifs (respectivement {\it espace des courbes non paramÈtrÈes} et {\it espace des courbes non paramÈtrÈes non orientÈes}), $\pi_O$ et $\pi_{SO}$ les projections associÈes et $\pi_S$ la projection naturelle de $\Pi S^2$ sur $\tilde{\Pi}S^2$. La fonctionnelle Ènergie $E_g$ initialement dÈfinie sur $\Lambda S^2$ par 
$$
E_g(c)=1/2 \int_{S^1} g_{c(t)}(\dot{c}(t),\dot{c}(t)) dt
$$ 
induit une application sur $\Pi S^2$ et $\tilde{\Pi}S^2$ respectivement, pour laquelle les points critiques sont exactement les orbites de points critiques de $\Lambda S^2$ pour l'action correspondante. Or les points critiques de $E$ dans $\Lambda S^2$ sont exactement les gÈodÈsiques fermÈes (Èventuellement rÈduites ‡ un point).

Soit $w$ une classe d'homologie non nulle de $\Lambda S^2$ (respectivement $\Pi S^2$ et $\tilde{\Pi}S^2$) ‡ coefficients quelconques. Pour chaque cycle singulier $u \in w$, on note $|u|$ la rÈunion des images des simplexes singuliers de $u$.

Nous avons alors le rÈsultat suivant :

\begin{theorem}(voir \cite{KLIN}) \label{klin}
Soit $\kappa_w=\inf_{u \in w} \sup_{c \in |u|} E_g(c)$. Alors $\kappa_w>0$ et il existe une gÈodÈsique fermÈe $c_w \in \Lambda M$ vÈrifiant $E_g(c_w)=\kappa_w$.
\end{theorem}

\bigskip

\noindent {\bf Familles particuliËres de courbes sur $(S^2,g_0)$.} La gÈomÈtrie de la sphËre canonique prÈsente naturellement plusieurs familles de courbes intÈressantes, dans le sens o˘ ces familles gÈnËrent des classes homologiques non nulles dans un des espaces de courbes dÈcrits dans le paragraphe prÈcÈdent. Nous considÈrons dans ce paragraphe $(S^2,g_0)$ comme la sphËre unitÈ de $\R^3$ pour le produit scalaire canonique, munie de la mÈtrique induite.

Notons tout d'abord $A S^2$ l'{\it espace des cercles paramÈtrÈs}. Un tel cercle est soit une application constante, soit un plongement $c :S^1 \rightarrow S^2$ paramÈtrÈ proportionnellement ‡ la longueur d'arc, dont l'image est l'intersection de $S^2$ avec un plan ‡ distance $<1$ de l'origine de $\R^3$. L'{\it espace des grands cercles} est notÈ $B S^2$.

L'image $\pi_O(B S^2)$ (espace des grands cercles non paramÈtrÈs) est clairement en bijection avec $S^2$, et nous pouvons dÈcrire la bijection comme suit : ‡ chaque $u \in S^2$ correspond le grand cercle non paramÈtrÈ de $S^2$ notÈ $\gamma(u)$ obtenu comme l'intersection du plan orthogonal ‡ $u$ avec $S^2$ orientÈ par la rËgle du tire-bouchon (cf figure 1).

\begin{figure}[ht]
\begin{center} 
\includegraphics[width=5cm]{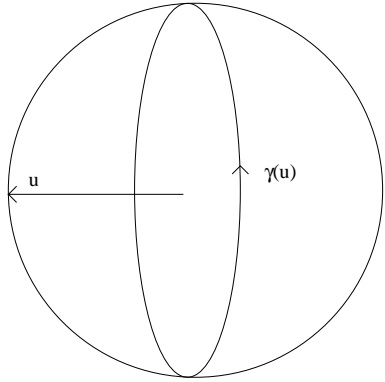}
\caption{Le grand cercle associÈ ‡ $u$}
\end{center}
\end{figure}

\medskip

Soit $\F=\{\pi_S(\gamma(\cos(t\pi),\sin(t\pi),0)), t\in [0,1]\}$ la famille de courbes paramÈtrÈe par l'intervalle $I=[0,1]$. Cette famille fournit un ÈlÈment $[F]\neq 0$ de l'homologie $H_1(\tilde{\Pi} S^2,\Z_2)$.

\medskip

Soit $u \in S^2$. Pour chaque $s \in [-1,1]$, nous pouvons associer un ÈlÈment $\gamma(u,s)$ de $\pi_O(A S^2)$ (espace des cercles non paramÈtrÈs) obtenu comme l'intersection avec $S^2$ du plan orthogonal ‡ $u$ translatÈ du vecteur $su$ et orientÈe naturellement (cf figure 2).

\begin{figure}[h]
\begin{center}
\includegraphics[width=8cm]{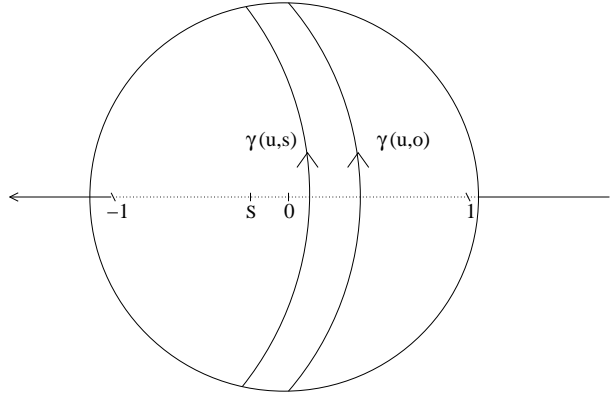}
\caption{La famille $\gamma(u,s)$}
\end{center}
\end{figure}

\medskip

Soit $\G_u$ la famille d'ÈlÈments de $\Pi S^2$ paramÈtrÈe par $S^1$ obtenu comme la concatÈnation de la famille $\{\gamma(u,s)) \mid s \in [-1,1]\}$ avec la famille de courbes rÈduites ‡ un point dÈfinie par un demi grand cercle reliant $\gamma(u,-1)$ ‡ $\gamma(u,1)$. Cette famille fournit un ÈlÈment non trivial de $H_1(\Pi S^2,\Z)$.

\bigskip

\noindent{\bf Tangent unitaire et mÈtrique associÈe.} Soit $U^{g_0} S^2$ le fibrÈ tangent unitaire ‡ $(S^2,g_0)$ et $\pi_{UM}$ la projection associÈe. Cet espace est isomorphe ‡ $BS^2$.

Clairement, on a une action libre de $S^1\simeq \R/(2\pi \Z)$ sur cet espace, dont le quotient coÔncide avec l'espace des grands cercles non paramÈtrÈs \linebreak $\pi_O(B S^2)\simeq S^2$.

La variÈtÈ $U^{g_0} S^2$ possËde une mÈtrique naturelle $g_1$ vÈrifiant 
$$
dv_{g_1}=\alpha_{g_0} \wedge \pi_O^*(dv_{g_0}), \eqno (2.3)
$$
o˘ $\alpha_{g_0}$ dÈsigne le tirÈ en arriËre de la $1$-forme de Liouville par l'isomorphisme musical bÈmol associÈ ‡ $g_0$, ainsi que l'ÈgalitÈ suivante pour toute fonction \linebreak $f \in C^1(U^{g_0} S^2,\R)$ :
$$
\int f dv_{g_1}=\int f d\sigma \times dv_{g_0}, \eqno (2.4)
$$
o˘ $\sigma$ dÈsigne la mesure canonique sur le cercle unitÈ $S^1$ (voir \cite{BESS}).

\subsection{DÈmonstration de l'inÈgalitÈ (2.2)}

Dans ce qui suit, et Ètant fixÈs $u \in S^2$ et $s \in [-1,1]$, nous allons considÈrer les longueurs des courbes $\gamma(u)$ et $\gamma(u,s)$, qui seront alors vues comme des ÈlÈments de $BS^2$ et $AS^2$ respectivement, et donc paramÈtrÈes proportionnellement ‡ la longueur d'arc.

\bigskip

Posons, pour tout $(q,v) \in U^{g_0} M$,
$$
\Lambda(q,v)=\sqrt{g_t(q)(v,v)}=1+tf(q).
$$

En intÈgrant cette fonction sur $U^{g_0} S^2$ par rapport ‡ la mesure canonique associÈe ‡ $g_1$, on obtient d'une part,
\begin{eqnarray} 
\nonumber   \int_{U^{g_0} S^2} (1+tf) dv_{g_1} & = & \int_{S^2} \left( \int_{S^1} (1+tf) d\sigma \right)dv_{g_0} \text{ (par (2.4)) }\\
\nonumber & = & 2\pi \int_{S^2} (1+tf) dv_{g_0}\\
\nonumber &=& 8 \pi^2, \\
\nonumber
\end{eqnarray}
et d'autre part, 
\begin{eqnarray}
\nonumber \int_{U^g S^2} (1+tf) dv_{g_1} & = & \int_{\pi_O(B S^2)} \left( \int_{\pi_O^{-1}(\gamma)} (1+tf) \cdot \alpha_g \right)dv_{g_0} \text{ (par (2.3)) }\\
\nonumber    & = & \int_{S^2} l_{g_t}(\gamma(u)) dv_{g_0}. \\
\nonumber 
\end{eqnarray} 

On en dÈduit donc que la moyenne de la longueur respectivement ‡ la mÈtrique $g_t$ sur l'espace des grands cercles non paramÈtrÈs est constante, ‡ savoir 
$$
\frac{1}{4\pi} \int_{S^2} l_{g_t}(\gamma(u)) dv_{g_0} =2\pi.
$$

\bigskip

Remarquons qu'un calcul simple nous montre que pour tout $u \in S^2$,
$$
l_{g_t}(\gamma(u))=2\pi + t \int_{\gamma(u)} f \cdot \alpha_{g_0}.
$$
Deux cas se prÈsentent donc ‡ nous :

\medskip

\noindent{\it Premier cas.} Pour tout $u \in S^2$, 
$$
\int_{\gamma(u)} f \cdot \alpha_{g_0}=0.
$$
 Alors pour tout $u \in S^2$ et $|t|<a$, $l_{g_t}(\gamma(u))=2\pi$. Comme toute courbe \linebreak $c : S^1=[0,1]/\{0,1\} \rightarrow S^2$ de $\Lambda S^2$ vÈrifie 
$$
2\cdot E_g(c)\leq l_g(c),
$$
 on obtient que pour tout $u \in S^2$ et tout $t \in ]-a,a[$, $E_{g_t}(\gamma(u)) \leq \pi$. Donc, pour tout $c \in \F$, 
$$
E_{g_t}(c) \leq \pi.
$$
 Comme $[\F] \neq 0$, on obtient par le thÈorËme \ref{klin} l'existence d'une gÈodÈsique fermÈe correspondant ‡ cette classe d'Ènergie au plus $\pi$, donc de longueur au plus $2\pi$ pour la mÈtrique $g_t$ avec $|t|<a$.
 
D'o˘ l'inÈgalitÈ (2.2) dans ce cas.

\medskip

\noindent{\it Second cas.} Il existe $u_0 \in S^2$, tel que 
$$
\int_{\gamma(u_0)} f \cdot \alpha_{g_0}<0.
$$
Remarquons qu'alors, il existe nÈcessairement $u_1 \in S^2$ tel que
$$
\int_{\gamma(u_1)} f \cdot \alpha_{g_0}>0.
$$

Pour tout $t \in ]0,a[$, $l_{g_t}(\gamma(u_0)) < 2\pi$. Donc il existe $\beta>0$ tel que pour tout $s \in ]-\beta,\beta[$ et tout $t \in ]0,a[$, 
$$
l_{g_t}(\gamma(u_0,s)) < 2\pi.
$$

 On se fixe alors $\alpha_0 \in ]0,a[$ tel que pour $t \in ]0, \alpha_0[$, $l_{g_t}(\gamma(u_0,s)) < 2\pi$ pour tout $s\in [-1,1]\setminus ]-\beta,\beta[$. Nous obtenons ainsi, pour $t \in ]0,\alpha_0[$ et pour tout $s\in [-1,1]$, 
$$
l_{g_t}(\gamma(u_0,s)) < 2\pi.
$$
 Donc pour chaque $c \in \G_{u_0}$, $E_{g_t}(c)<\pi$. De la non trivialitÈ de $[\G_{u_0}]$, nous obtenons par le thÈorËme \ref{klin} l'existence d'une gÈodÈsique fermÈe correspondant ‡ la classe $[\G_{u_0}]$ d'Ènergie au plus $\pi$, donc de longueur au plus $2\pi$ pour la mÈtrique $g_t$ et ce pour tout $t \in ]0,\alpha_0[$.
 
 \medskip
 
 De maniËre analogue, nous obtenons l'existence d'un $\alpha_1 \in ]0,a[$ tel que pour tout $t \in ]-\alpha_1,0[$, il existe une gÈodÈsique fermÈe correspondant ‡ la classe $[\G_{u_1}]$ d'Ènergie au plus $\pi$, donc de longueur au plus $2\pi$ pour la mÈtrique $g_t$.

\medskip

On pose alors $\alpha=\min\{\alpha_0,\alpha_1\}$ et l'inÈgalitÈ (2.2) est dÈmontrÈe pour $t \in ]-\alpha,\alpha[$.

\section{Aire systolique et paritÈ}

Nous prÈsentons ici deux rÈsultats : le premier est d˚ ‡ V. Guillemin \cite{guil76} et le second ‡ P. Pu \cite{pu52}. Ces rÈsultats nous permettent de complÈter le panorama des connaissances actuelles sur le comportement local de l'aire systolique au voisinage de $g_0$, ce qui justifie leur prÈsentation ici.

\medskip

Plus prÈcisÈment, notons $a$ l'antipodie de $(S^2,g_0)$. On dÈfinit le {\it sous-espace des fonctions paires} 
$$
C^\infty_+(S^2,\R)=\{f \in C^\infty(S^2,\R) \mid f \circ  a=f\},
$$
et le {\it sous-espace des fonctions impaires}
$$
C^\infty_-(S^2,\R)=\{f \in C^\infty(S^2,\R) \mid f \circ a=-f\}.
$$

Nous avons la dÈcomposition suivante :
$$
C^\infty(S^2,\R)=C^\infty_+(S^2,\R) \oplus C^\infty_-(S^2,\R).
$$

\bigskip

Nous allons voir que :

\medskip

\noindent 1). Pour toute fonction $f \in C^\infty_-(S^2,\R)$, il existe une variation lisse $g_t$ de $g_0$ par des mÈtriques lisses dans la direction tangente prescrite par $f$ telle que pour $t\geq 0$,
$$
\S(S^2,g_t)=1/\pi.
$$

\noindent 2). Pour toute mÈtrique lisse $g$ paire respectivement ‡ l'antipodie ({\it i.e.} vÈrifiant $a^\ast g =g$), nous avons
$$
\S(S^2,g)\geq 1/\pi,
$$
avec ÈgalitÈ si et seulement si $g=g_0$.

Donc pour toute fonction $f$ non nulle de $C^\infty_+(S^2,\R)$, la variation \linebreak  $g_t=(1+tf)\cdot g_0$ de $g_0$ dans la direction tangente prescrite par $f$ vÈrifie  pour $t\neq 0$,
$$
\S(S^2,g_t)>1/\pi.
$$

\bigskip

Remarquons qu'en combinant les points 1). et 2)., nous retrouvons le corollaire prÈsentÈ dans l'introduction.

\bigskip

\noindent {\bf Fonctions impaires et mÈtriques de Zoll.} La dÈfinition suivante est classique (voir \cite{BESS}). Une mÈtrique riemannienne $g$ de $S^2$ est dite une {\it mÈtrique de Zoll} si toutes les g\'eod\'esiques de $(S^2,g)$ sont p\'eriodiques de p\'eriode $2\pi$.  Remarquons que dans ce cas, $\s(S^2,g)=2\pi$.
Il est clair que la mÈtrique standard $g_0$ vÈrifie cette propriÈtÈ : les gÈodÈsiques fermÈes sont les ÈlÈments de $BS^2$. Une mÈtrique de Zoll $g$ vÈrifie l'ÈgalitÈ suivante (voir \cite{BESS} et \cite{wein74}) :
$$
\A(S^2,g)=\A(S^2,g_0),
$$
ce qui a pour consÈquence que
$$
\S(S^2,g)=\frac{1}{\pi}.
$$

A l'heure actuelle, les exemples connus de mÈtriques de Zoll sur la $2$-sphËre, autres que la mÈtrique canonique, sont les suivants :

\noindent - les mÈtriques de Zoll de rÈvolution (voir \cite{BESS} et \cite{zoll03}),

\noindent - certaines mÈtriques obtenues comme variation lisse de la mÈtrique standard $(S^2, g_0)$ :

\begin{theorem} (voir \cite{funk13} et \cite{guil76}) \label{funk}
 Pour toute fonction $f \in C^\infty_-(S^2,\R)$, il existe une famille lisse ‡ un paramËtre de $C^\infty$-fonctions $\rho_t$ telle que $\rho_0=0$, $\frac{d \rho_t}{dt}\mid_0=f$ et $\exp{(\rho_t)}\cdot g_0$ soit une mÈtrique de Zoll.
\end{theorem}

Nous obtenons donc que pour toute fonction impaire $f \in C^\infty_-(S^2,\R)$, il existe une variation $\Psi_t^- \cdot g_0$ de la mÈtrique $g_0$ telle que $\frac{d \Psi_t^-}{dt}\mid_0=f$ et
$$
\S(S^2,\Psi_t^- \cdot g_0)=\frac{1}{\pi}.
$$
Ceci dÈmontre le point 1).

\bigskip

\noindent {\bf Fonctions paires et thÈorËme de Pu.} Soit $g$ une mÈtrique sur $S^2$ lisse et paire respectivement ‡ l'antipodie. Comme $a^\ast g =g$, nous pouvons former le quotient $S^2/a\simeq\R P^2$ et le munir de la mÈtrique induite par $g$ que nous noterons $\tilde{g}$. 

…tant donnÈe une mÈtrique riemannienne $\hat{g}$ sur $\R P^2$, nous dÈfinissons la systole de $(\R P^2,\hat{g})$ comme la plus petite longueur d'une courbe non contractile et la notons $\s(\R P^2,\hat{g})$. Cette longueur est rÈalisÈe comme la longueur d'une gÈodÈsique fermÈe non contractile. Nous pouvons alors dÈfinir l'aire systolique de $(\R P^2,\hat{g})$ comme le quotient
$$
\S(\R P^2,\hat{g})=\frac{\A(\R P^2,\hat{g})}{(\s(\R P^2,\hat{g}))^2}.
$$

Nous observons que pour la mÈtrique $\g$ dÈfinie ‡ partir de la mÈtrique $g$ ci-dessus, nous avons
$$
\A(\R P^2,\g)=\frac{1}{2}\A(S^2,g),
$$
et 
$$
\s(\R P^2,\g)\geq \frac{1}{2}\s(S^2,g).
$$

On en dÈduit donc :
$$
\S(S^2,g)\geq \frac{\S(\R P^2,\g)}{2}.
$$

L'aire systolique du plan projectif vÈrifie l'inÈgalitÈ suivante :

\bigskip

\begin{theorem} (voir \cite{pu52})
 Soit $\hat{g}$ une mÈtrique riemannienne sur $\R P^2$. Alors
$$
\S(\R P^2,\hat{g})\geq \frac{2}{\pi},
$$
avec ÈgalitÈ si et seulement si $\hat{g}$ est la mÈtrique standard.
\end{theorem}

Nous en dÈduisons donc pour toute mÈtrique paire $g$, 
$$
\S(S^2,g)\geq \S(S^2,g_0)=\frac{1}{\pi},
$$
avec ÈgalitÈ si et seulement si $g=g_0$. 

Ceci dÈmontre le point 2).

\section{Remarques et perspectives}

Il serait intÈressant d'Èclairer plusieurs points :

\medskip

\noindent 1). Pour tout vecteur tangent non nul $f \in C^\infty(S^2,\R)$, nous avons construit une variation de mÈtrique $\Phi_t \cdot g_0$ telle que $\Phi_0=1$, $\left. \frac{d \Phi_t}{dt}\right |_0=f$, et pour laquelle \linebreak $\S(S^2,g_t) \geq \S(S^2,g_0)$ pour $t \geq 0$. Il semble raisonnable de penser que la mÈtrique $g_0$ est un minimum local de l'aire systolique dans l'espace des mÈtriques conformes ‡ $g_0$, mais il n'y a, ‡ l'heure actuelle et ‡ la connaissance de l'auteur, aucune dÈmonstration de cette conjecture.

\medskip

\noindent 2). De maniËre plus ambitieuse, il serait intÈressant de mieux comprendre le comportement local de l'aire systolique au voisinage de $g_0$. La troisiËme section suggËre une conjecture de la forme : il existe un voisinage ouvert $U_0$ de $0$ dans $ C^\infty(S^2,\R)=C^\infty_+(S^2,\R) \oplus C^\infty_-(S^2,\R)$, un voisinage $V_{g_0}$ de $g_0$ dans l'espace des mÈtriques lisses conformes ‡ $g_0$ muni de la topologie forte et un homÈomorphisme local $\phi : U_0 \rightarrow V_{g_0}$ tel que pour toute mÈtrique $g \in V_{g_0}$, si on note $\pi_+$ et $\pi_-$ les projections respectives de $C^\infty(S^2,\R)$ sur $C^\infty_+(S^2,\R)$ et $C^\infty_-(S^2,\R)$, nous ayons 
$$
\S(S^2, g)= \S(S^2,g_0)(1+||\pi_+(\phi^{-1}(g))||^2_2),
$$
o˘ $||f||_2=\sqrt{\int_{S^2} f^2 dv_{g_0}}$ pour tout $f \in C^\infty(S^2,\R)$.

\medskip

\noindent 3). Une question intÈressante est la suivante : toute mÈtrique de Zoll est-elle un point critique de la fonctionnelle aire systolique ? Plus gÈnÈralement, la dÈcouverte d'autres points critiques serait remarquable, et il n'est plas clair qu'elles sont les mÈtriques (en-dehors des mÈtriques de Zoll et de l'exemple de Calabi) qui fourniraient de bons candidats.

\bigskip
\bigskip

\noindent {\it Remerciements.} - L'auteur tient ‡ exprimer ici ces remerciements ‡ I. Babenko, pour lui avoir prÈsentÈ ce sujet de recherche, ainsi que pour les nombreuses conversations qui s'en sont suivies. Ses remerciements s'adressent Ègalement ‡ M. Berger pour ses rÈponses sur certains points, ‡ S. Sabourau pour ses remarques pertinentes, ainsi qu'‡ l'UniversitÈ de GenËve en la personne de P. de la Harpe d'une part, et au Centre Interfacultaire Bernouilli d'autre part, pour leur accueil respectif durant l'Èlaboration de cet article.

\end{document}